\documentstyle[amssymb,12pt]{article}

\begin{document}

\begin{titlepage}
\title{Some Ricci Flat (pseudo-) Riemannian Geometries}
\vspace{5cm}
\author{Metin G{\" u}rses \\
{\small Department of Mathematics, Faculty of Science} \\
{\small Bilkent University, 06533 Ankara - Turkey}\\
email:gurses@fen.bilkent.edu.tr}

\maketitle
\begin{abstract}
We define a class of two dimensional surfaces conformally related to minimal 
surfaces in  flat three dimensional
geometries. By the utility of the metrics of such surfaces we give 
a construction of the metrics of $2\,N$ dimensional Ricci flat (pseudo-) 
Riemannian geometries.

\end{abstract}
\end{titlepage}

\section{Introduction}

Let $(S, g)$ denote a two dimensional geometry where $S$ is a surface 
in a three dimensional flat manifold $M_{3}$ and $g$  is a (pseudo-) 
Riemannian metric on $S$ with a non vanishing determinant, $det(g)$.
Furthermore we assume that $g$ satisfies the following conditions

\begin{eqnarray}
\partial_{\mu}\, (g^{\mu \nu}\, g^{-1} \partial_{\nu}\,g)=0,\\
R+{1 \over 4}\, tr [g^{\mu \nu} \partial_{\mu}\, g^{-1} \partial_{\nu}\,g]=0.
\end{eqnarray}

\noindent
where $R$ is the Ricci scalar (Gaussian curvature) of $S$ (please see the
next section for our conventions). We shall see in the following sections that
some surfaces which are conformally related to minimal surfaces  satisfy
the above conditions.

The importance of such surfaces arises when we are interested in
even dimensional Ricci flat geometries. By the utility the metric $g$ of
these surfaces we shall give a  construction (without solving any  further
differential equations) of the  metric of a $2N$ dimensional 
Ricci flat (pseudo-) Riemannian geometries.

Ricci flat geometries are
important not only in differential geometry and general relativity
but also in gravitational instantons and in brane solutions of 
string theory \cite{brec}.

\section{Conformally related minimal surfaces}

Let $\phi$ be a differentiable function of $x^{1}$ and $x^{2}$ and $S_{0}$
be the surface in a three dimensional manifold $M_{3}$ (not necessarily
Euclidean) with a pseudo-Riemannian metric $g_{3}$ defined through
$ds^2=g_{0\, \mu \nu}\,dx^{\mu}\, dx^{\nu}+\epsilon \, (dx^{3})^{2}$, where
$\mu,\nu=1,2$ , $\epsilon= \pm 1$ and $g_{0}$ is a constant, invertible,
symmetric $2 \times 2$ matrix. In this work we assume Einstein summation
convention, i.e., the repeated indices are summed up. $S_{0}$ is given as 
the graph of the function $\phi$, i.e., $S_{0}=\{ (x^{1},x^{2},x^{3}) 
\in M_{3} | x^{3}= \phi (x^{1},x^{2}) \}$. Then the metric on $S_{0}$ 
is given by

\begin{equation}
h_{\mu \nu}=g_{0\, \mu \nu}+\epsilon\, \phi_{, \mu}\, \phi_{, \nu}.
\label{2met}
\end{equation}

\noindent
Since\, $det\, h = (det \, g_{0})\, \rho$ where

\begin{equation}
\rho=1+\epsilon\, g_{0}^{\mu \nu}\, \phi_{, \mu}\, \phi_{,\nu}
\end{equation}

\noindent
then $h$ is everywhere (except at those points where $\rho=0$) 
invertible. Its inverse is given by

\begin{equation}
h^{\mu \nu}=g_{0}\,^{\mu \nu}-{\epsilon \over \rho}\, \phi^{\mu}\, \phi^{\nu}
\end{equation}

\noindent
where $g_{0}^{\mu \nu}$ are the components of the inverse matrix $g_{0}^{-1}$
of $g_{0}$. Here the indices are lowered and raised by the metric $g_{0}$
and its inverse $g_{0}\, ^{-1}$ respectively. For instance , $\phi^{\mu}_{\nu}=
g_{0}^{\mu \alpha}\, \phi_{, \alpha \nu}$.

The Ricci tensor corresponding to the metric in  (\ref{2met}) is given by

\begin{equation}
r_{\mu \nu}={\epsilon \over \rho}\, (\nabla^{2}\, \phi)\, \phi_{,\mu \nu}
-{\epsilon \over \rho}\, \phi_{,\mu}\,^{\alpha}\, \phi_{,\nu \alpha}+
{1 \over 4 \rho^{2}}\, \rho_{,\mu}\, \rho_{, \nu}, \label{2ric}
\end{equation}

\noindent
where

\begin{equation}
\nabla^{2}\, \phi= h^{\mu \nu}\, \phi_{, \mu \nu}=
g_{0}\,^{\mu \nu}\, \phi_{, \mu \nu}-{ 1 \over  2 \rho}\, \phi^{\alpha}\,
\rho_{, \alpha}
\end{equation}

\noindent
The Ricci scalar or the Gaussian curvature $K$ and the minimal curvature
$H$ are obtained as

\begin{eqnarray}
K&=&{\epsilon \over \rho^2}\, [(\phi^{\alpha}\,_{\alpha})^2-\phi^{\alpha \beta}\,
\phi_{,\alpha \beta}]\\
H&=&{1 \over \sqrt{\rho}}\, h^{\mu \nu}\, \phi_{,\mu \nu}, \label{min}
\end{eqnarray}

\noindent
The following proposition is valid only for the case of two dimensional
geometries.

\noindent
{\it Proposition 1.\,
$\phi_{,\alpha \mu}\, \phi_{,\beta \gamma}-\phi_{,\alpha \beta}\,
\phi_{,\mu \gamma}=
-\lambda_{0}\,(g_{0\, \alpha \mu}\,g_{0\,\beta \gamma}-
g_{0\,\alpha \beta}\,g_{0\,\gamma \mu})$

\noindent
where

\begin{equation}
\lambda_{0}={1 \over 2}\,[\phi^{\alpha \beta}\, \phi_{\alpha \beta}-
(\phi^{\alpha}_{\alpha})^2].
\end{equation}
}
\noindent
The following corollaries will be very useful in this section

\noindent
{\it Corollary 1.\,
$\phi^{\alpha}_{\mu}\, \phi_{,\alpha \nu}-
\phi^{\alpha}_{\alpha}\, \phi_{,\mu \nu}=\lambda_{0}\,g_{0\,\mu \nu}$
}\\
\noindent
{\it Corollary 2.
$r_{\alpha \beta}={K \over 2}\, h_{\alpha \beta},~~~\lambda_{0}=
-{\epsilon \over 2}\,\rho^2\,K$
}\\
\noindent
For the minimal surfaces we have $H=0$ and the following result

\noindent
{\it Proposition 2.\, If $H=0$ then (\cite{gur1})
\begin{eqnarray}
\partial_{\alpha}\,[\sqrt{\rho}\, h^{\alpha \beta}\,
\partial_{\beta}\, \phi]&=&0,\\
\partial_{\alpha}\,( \sqrt{\rho}\, h^{\alpha \beta})&=&0.
\end{eqnarray}
}
\noindent
We now define surfaces which are locally conformal 
to minimal surfaces. Let $S$ be such a surface, i.e., locally conformal
to $S_{0}$. Then the  metric on $S$ is given by

\begin{equation}
g_{\alpha \beta}={ 1 \over \sqrt{\rho}}\, h_{\alpha \beta}.
\end{equation}

\noindent
It is clear that $det\, g= det \,g_{0} \ne 0$. In the sequel we shall assume
that the surface $S_{0}$ is minimal and hence the metric defined on it
satisfies all the equivalent conditions in proposition 2. The corresponding
Ricci tensor of $g$ is given as

\begin{equation}
R_{\alpha \beta}=r_{\alpha \beta}-(\nabla_{g}^{2}\, {\psi_{0}})\,
g_{\alpha \beta},
\end{equation}

\noindent
where $\psi_{0}=-{1 \over 4}\, log (\rho)$ and $\nabla_{g}^{2}$ is 
the Laplace-Beltrami operator with respect to
the metric $g$. Using the above results we have

\noindent
{\it Proposition 3.\, The following are some identities related
to the conformally related surface $S$.
\begin{eqnarray}
R&=&-{1 \over 4} g^{\alpha \beta}\, tr [\partial_{\alpha}\, g^{-1}\, 
\partial_{\beta}\, g]\\
R_{\alpha \beta}&=&-{\rho-1 \over 2}\, r_{\alpha \beta},\\
R&=&\sqrt{\rho}\,r-2 \nabla_{h}^{2}\, \psi_{0}.
\end{eqnarray}

\noindent
Here $g$ is the $2 \times 2$ matrix of $g_{\alpha \beta}$ and $g^{-1}$
is its inverse. The operation tr is the standard trace operation for
matrices.}

\noindent
Let $v_{\alpha}=(1,0),~v^{\prime}_{\alpha}=(0,1)$ and
$u^{\alpha}=(1,0),~u^{\prime\, \alpha}=(0,1)$. We now define some functions
over $S$.

\begin{eqnarray}
\xi_{1}&=&g^{\alpha \beta}\, v_{\alpha}\, v_{\beta},~~~\xi_{2}=
g^{\alpha \beta}\, v^{\prime}_{\alpha}\, v^{\prime}_{\beta},\\
w_{1}&=&\sqrt{\rho}\, g_{\alpha \beta}\,u^{\alpha}\, u^{\beta},~~~
w_{2}=\sqrt{\rho}\, g_{\alpha \beta}\, u^{\prime \, \alpha}
\,u^{\prime\, \beta}.
\end{eqnarray}

\noindent
It is now easy to prove

\noindent
{\it Proposition 4.
\begin{eqnarray}
\nabla_{g}\,^2 \zeta -a_{0}\, R&=&-a_{0}\, \sqrt{\rho}\,K ,
\label{oz01}\\
\nabla_{g}\,^2 \psi_{1}-(a_{1}+a_{2})\, R&=&0\\
\nabla_{g}\,^2 \psi_{2}-2(b_{1}+b_{2})\, R&=&-(b_{1}+b_{2})\,\sqrt{\rho}\,K,
\end{eqnarray}

\noindent
where

\begin{eqnarray}
\zeta&=&{a_{0} \over 2}\, log (\rho),\\
\psi_{1}&=&a_{1}\,log(\xi_{1})+a_{2}\, log(\xi_{2}),\\
\psi_{2}&=&b_{1}\, log(w_{1})+b_{2}\, log(w_{2}).
\end{eqnarray}

\noindent
Here $a_{0},a_{1},a_{2},b_{1},$ and $b_{2}$ are arbitrary constants.}

\noindent

There is another function $\mu=(b_{1}+b_{2})\, \zeta-a_{0}\, \psi_{2}$
satisfies similar equation as $\psi_{1}$

\begin{equation}
\nabla_{g}\,^{2}\, \mu=-a_{0}\,(b_{1}+b_{2})\,R, \label{mu}
\end{equation}

\noindent
Using Eq.(\ref{oz01}) the functions $\mu$ and $\psi_{1}$ satisfy
a similar type of equation

\begin{equation}
\nabla_{g}^{2} \sigma=-{c \over 4}\, g^{\alpha \beta}\, tr [
(\partial_{\alpha}\,g^{-1})\, \partial_{\beta}\, g], \label{cc1}
\end{equation}

\noindent
where $c=(a_{1}+a_{2})$ when $\sigma=\psi_{1}$ and 
$c=-a_{0}(b_{1}+b_{2})$ when $\sigma=\mu$.

It is easy to show that 

\begin{equation}
\xi_{1}={w_{2} \over det\, g_{0}\, \sqrt{\rho}},~~ 
\xi_{2}={w_{1} \over det\, g_{0}\,\sqrt{\rho}}
\end{equation}

\noindent
Hence $\psi_{1}$ will not be considered as an independent function.
It is  interesting and important to note that under the minimality
condition the matrix $g$ satisfies the following condition as well.

\noindent
{\it Proposition 5.\, Minimality of $S_{0}$, $H=0$, also implies
a sigma model \cite{gur1}, \cite{gur2} like equation for $g$, i.e.,

\begin{equation}
\partial_{\alpha}\,[g^{\alpha \beta}\, g^{-1}\, \partial_{\beta}\,g]=0
\end{equation}
}
\section{Four Dimensions}

Let the metric of a four dimensional manifold $M_{4}$ be given by

\begin{equation}
ds^2=e^{2\psi}\,g_{\alpha \beta}\,dx^{\alpha}\,dx^{\beta}+
\epsilon_{1}\,g_{\alpha \beta}\,dy^{\alpha}\,dy^{\beta}, \label{4met}
\end{equation}

\noindent
where $\psi$ is a function of $x^{\alpha}$ and $\epsilon_{1}=\pm 1$.
Local coordinate of $M_{4}$ are denoted as $x^{a}=(x^{\alpha},
y^{\alpha}), ~a=1-4$

\noindent
{\it Proposition 6.\, The Ricci flat equations $R_{ab}=0$ for the metric
(\ref{4met}) are given in two sets. One set satisfied identically
due to the Proposition 5 above and the second one is given by

\begin{equation}
\nabla_{g}\,^{2} \psi=0
\end{equation}
}

\noindent
There are two independent functions satisfying the above
Laplace equation , $\phi$ and $\mu$.
Using (\ref{mu}) we find that $\psi=e_{0}\, \phi+e_{1}\, \mu$
where $e_{0}$ and $e_{2}$ are arbitrary constants and $b_{2}=-b_{1}$.
Combining all these constants we find that

\begin{equation}
e^{2\psi}=e^{2 e_{0} \, \phi}\,w_{1}^{-2m_{1}}\,w_{2}^{-2m_{2}},
\label{psi}
\end{equation}

\noindent
where $m_{1}$ and $m_{2}$ are constants satisfying $m_{1}+m_{2}=0$.
Then the line element (\ref{4met}) becomes

\begin{equation}
ds^{2}={e^{2\,e_{0}\, \phi} \over w_{1}^{2m_{1}}\,w_{2}^{2m_{2}}} 
\,{h_{\alpha \beta}\,dx^{\alpha}\,dx^{\beta} \over \sqrt{\rho}}+
{h_{\alpha \beta}\,dy^{\alpha}\,dy^{\beta} \over \sqrt{\rho}}, \label{ins}
\end{equation}

\noindent
where $\phi$ satisfies the minimality condition ($H=0$)  (\ref{min})  which 
is explicitly given by

\begin{equation}
[k_{2}+\epsilon\,(\phi_{,y})^{2}]\,\phi_{,xx}-2[k_{0}+
\epsilon\, \phi_{,x}\, \phi_{,y}]\, \phi_{,xy}+[k_{1}+\epsilon\,
(\phi_{,x})^{2}]\, \phi_{,yy}=0, \label{min1}
\end{equation}

\noindent
where we take $(g_{0})_{ 11}=k_{1},~(g_{0})_{01}=k_{0},~(g_{0})_{22}=k_{2}$
and   assume that $det\, (g_{0})=k_{1}\,k_{2}-k_{0}^{2} \ne 0$. Hence
the functions $w_{1}$ and $w_{2}$ are given explicitly as

\begin{equation}
w_{1}=k_{1}+\epsilon\, (\phi_{,x})^{2},~~w_{2}=k_{2}+
\epsilon\, (\phi_{,y})^{2}. \label{w12}
\end{equation}

\noindent
The metric in (\ref{ins}) with  $e_{0}=0, m_{1}=m_{2}=0$ reduces to 
an instanton metric \cite{yn}. 

\section{Higher Dimensions}

Let $M_{2+2n}$ be a $2+2n$ dimensional manifold with a metric

\begin{equation}
ds^2=e^{2\Phi}\,g_{\alpha \beta}\,dx^{\alpha}\,dx^{\beta}
+G_{AB}\,dy^{A}\,dy^{B}, \label{nmet}
\end{equation}

\noindent
where the local coordinates of $M_{2+2n}$ are given by $x^{\alpha+A}=
(x^{\alpha},y^{A}),~A=1,2, \cdots ,2n$, $\Phi$ and $G_{AB}$ are functions
of $x^{\alpha}$ alone. The Einstein  equations are given in
the following following proposition

\noindent
{\it Proposition 7.\, The Ricci flat equations for the metric in
(\ref{nmet}) are given by

\begin{eqnarray}
\partial_{\alpha}\,[g^{\alpha \beta}\,G^{-1}\, \partial_{\beta}\,G]=0,\\
\nabla_{g}\,^{2}\, \Phi={1 \over 8}\, g^{\alpha \beta}\,
tr [(\partial _{\alpha} G^{-1}) \, \partial_{\beta}\,G]+{R \over 2},
\end{eqnarray}

\noindent
where $G$ is $2n \times 2n$ matrix of $G_{AB}$ and $G^{-1}$ is
its inverse.}

\noindent
Let us choose $G$ as a block diagonal matrix and
each block is the $2 \times 2$ matrix $g$. This means that
the metric in (\ref{nmet}) reduces to a special form

\begin{eqnarray}
ds^2=e^{2\Phi}\, g_{\alpha \beta}\,dx^{\alpha}\,dx^{\beta}+
\epsilon_{1}\,g_{\alpha \beta}\,dy_{1}^{\alpha}\,dy_{1}^{\beta}+
\cdots +\epsilon_{n}\,g_{\alpha \beta}\,dy_{n}^{\alpha}\,
dy_{n}^{\beta}, \label{smet}
\end{eqnarray}

\noindent
where the local coordinates of $M_{2+2n}$ are given by $x^{\alpha+A}=
(x^{\alpha},y_{1}^{\alpha}, \cdots , y_{n}^{\alpha})$, $\epsilon_{i}
=\pm 1,~ i=1,2, \cdots , n$. Then we have the following  theorem

\noindent
{\it Theorem.\, For every two dimensional minimal surface $S_{0}$
immersed in a three
dimensional manifold $M_{3}$ there corresponds a $2N=2+2n$-dimensional
Ricci flat (pseudo-) Riemannian geometry with the metric given in (\ref{smet})
with

\begin{equation}
e^{2\, \Phi}=e^{2\, \psi}\, \,{w_{1}^{-2n_{1}}\, w_{2}^{-2n_{2}}}
\,\, \rho^{n_{1}+n_{2} } , \label{phi}
\end{equation}

\noindent
where $\psi$ is given in (\ref{psi}), $w_{1}$ and $w_{2}$ are given
in (\ref{w12}), $n_{1}$ and $n_{2}$ satisfy

\begin{equation}
n_{1}+n_{2}={n-1 \over 2}. \label{n12}
\end{equation}
}

\noindent
{\it Proof of the theorem:}\, Using proposition 7 for the metric
(\ref{smet}) the Ricci flat equations reduce to the following
equation

\begin{equation}
\nabla_{g}\,^{2}\, \Phi={n-1 \over 8}\, g^{\alpha \beta}\, tr
[ (\partial_{\alpha} g^{-1})\, \partial_{\beta} g]
\end{equation}

\noindent
Hence, using (\ref{cc1}) and letting $a_{0}\,b_{1}=n_{1},~
a_{0}\,b_{2}=n_{2}$ and $\Phi=\mu$ we find (\ref{phi}) with the condition
(\ref{n12}).
All metric functions $\psi$, $w_{1}$, $w_{2}$ and $g_{\alpha \beta}$
are expressed explicitly in terms the function $\phi$ and its
derivatives $\phi_{,x}$ and $\phi_{,y}$. This means that for each solution
$\phi$ of (\ref{min1}) there exists a $2N$-dimensional metric (\ref{smet}).
This completes the proof of the theorem.

The dimension of the manifold is $2N=4+4(n_{1}+n_{2})$.
Here $n=1$ corresponds to the four dimensional case. 
The signature of the geometry depends on the 
signature of $S$. If the signature of $S$ is zero then the signature of 
$M_{2N}$ is also zero. If the signature of $S$ is $2$ then the signature 
of $M_{2N}$ is $2\,(1+\epsilon_{1}+ \cdots +\epsilon_{n})$.

\vspace{1cm}

This work is partially supported by the Scientific and Technical 
Research Council of Turkey (TUBITAK) and Turkish Academy of Sciences (TUBA).

\end{document}